\documentclass[11pt]{article}
\usepackage{amsfonts}
\usepackage{latexsym}
\usepackage{amsmath}
\usepackage{amssymb}
\usepackage{amsthm}
\newtheorem{theorem}{Theorem}[section]
\newtheorem{lemma}[theorem]{Lemma}

\newtheorem{proposition}[theorem]{Proposition}
\newtheorem{corollary}[theorem]{Corollary}

\theoremstyle{definition}
\newtheorem{definition}[theorem]{Definition}

\theoremstyle{remark}

\newtheorem*{note*}{Note}



\begin{document}

\title{\bf 
$f$-Divergence  for convex bodies }

\author{ Elisabeth M. Werner 
\thanks{Partially supported by an NSF grant, a FRG-NSF grant and  a BSF grant}}

\date{}

\maketitle

\begin{abstract}

We introduce $f$-divergence, a concept from  information theory and statistics, 
for  convex bodies in $\mathbb{R}^n$. We prove that $f$-divergences are $SL(n)$ invariant valuations
and we establish an affine isoperimetric inequality for these quantities.
We show that generalized affine surface area and in particular the $L_p$ affine surface area
from the $L_p$ Brunn Minkowski theory are special cases of $f$-divergences.
\end{abstract}

 \vskip 2cm

\section{Introduction.}

In information theory, probability theory and statistics, an $f$-divergence is a function $D_f(P,Q)$  that measures the difference between two probability distributions $P$ and $Q$. The divergence is intuitively an average, weighted by the function $f$, of the odds ratio given by $P$ and $Q$.
These divergences were introduced  independently by Csisz\'ar \cite{Csiszar}, Morimoto \cite{Morimoto1963}  and Ali \& Silvey \cite{AliSilvery1966}. Special cases of $f$-divergences are the Kullback Leibler divergence or relative entropy and the 
R\'enyi divergences (see Section 1). 
\par
Due to a number of highly influential works (see, e.g., 
\cite{Ga1} - \cite{HLYZ}, 
\cite{Klain1}, \cite{Klain2},
\cite{Lud2}, \cite{Lud3}, 
\cite{LR1} -  \cite{Lu1988},  
\cite{LYZ2002}, \cite{LYZ2004}, 
\cite{LZ} - \cite{MW2}, 
\cite{NPRZ}, \cite{RuZ}, 
\cite{Schu} - 
\cite{WernerYe2010},
\cite{Z1}  -  \cite{Z3}), the
$L_p$-Brunn-Minkowski theory is now  a central part of modern convex geometry.
A  fundamental notion within this theory 
is $L_p$ affine surface area,   introduced by Lutwak   
in the ground breaking paper \cite{Lutwak1996}. 
\par
It was shown in \cite{Werner2012} that $L_p$ affine surface areas  are  entropy powers of
R\'enyi divergences  of the cone measures of a convex body  and its polar, thus establishing 
further connections between information theory and convex geometric 
analysis. 
Further examples of such connections are e.g. several  papers by Lutwak, Yang, and Zhang \cite{LYZ2000, LYZ2002/1, LYZ2004/1, LYZ2005} and the recent article \cite{PaourisWerner2011} where it is shown how relative entropy appears in convex geometry. 
\par
In this paper we introduce $f$-divergences to the theory of convex bodies and thus 
 strengthen the already existing ties  between information theory and convex geometric 
analysis. 
We show that generalizations 
of the $L_p$ affine surface areas, the $L_\phi$ and $L_\psi$ affine surface areas introduced in 
\cite{LR2} and \cite{Ludwig2010},  are in fact $f$-divergences for special functions $f$.
We show that $f$-divergences are $SL(n)$ invariant valuations
and establish an affine isoperimetric inequality for these quantities. Finally, we give geometric characterizations
of $f$-divergences.

\vskip 4mm
Usually, in the literature, $f$-divergences are considered for convex functions $f$. A similar theory with the obvious modifications can be developed for concave functions. Here, we restrict ourselves to consider the convex setting.

\vskip 4mm
{\bf Further Notation.}
\vskip 2mm
\noindent
We work in ${\mathbb R}^n$, which is
equipped with a Euclidean structure $\langle\cdot ,\cdot\rangle $.
We write $B_2^n$
for the Euclidean unit ball centered at $0$ and $S^{n-1}$ for the
unit sphere. Volume is denoted by $|\cdot |$ or, if we want to emphasize the dimension,  by $\text{vol}_d(A)$ for a $d$-dimensional set $A$. 
\par
Let $\mathcal{K}_0$ be the space of convex bodies $K$ in $\mathbb{R}^n$ that contain the origin in their interiors. 
Throughout the paper, we will only consider  such $K$. For $K \in \mathcal{K}_0$, 
$K^\circ=\{ y\in \mathbb{R}^n: \langle x, y \rangle \leq 1 \  \text{for all }  \  x \in K\}$ is the polar body of $K$.
For a point $x \in \partial K$, the boundary of $K$,  $N_K(x)$ is the outer unit normal 
in $x$ to $K$ and $\kappa_K(x)$, or, in short $\kappa$, is the (generalized) Gauss curvature in $x$.
We  write $K\in
C^2_+$,  if $K $ has $C^2$ boundary $\partial K$  with everywhere strictly
positive Gaussian curvature $\kappa_K$.  By $\mu$ or $\mu_K$ we denote  the usual surface area measure on $\partial K$ and by
 $\sigma $  the usual surface area measure on $S^{n-1}$.
\par
Let  $K$ be a convex body in $\mathbb{R}^n$ and let  $u \in S^{n-1}$. Then
$h_{K}(u)$ is the support function of $K$ in direction $u\in S^{n-1}$,
and $f_{K}(u)$ is the curvature function, i.e. the reciprocal of
the Gaussian curvature $\kappa _K(x)$ at the point $x \in
\partial K$ that has $u$ as outer normal.

\vskip 5mm
\section{$f$-divergences.}

Let $(X, \mu)$ be a measure space  and let  $dP=pd\mu$ and  $dQ=qd\mu$ be probability  measures on $X$ that are  absolutely continuous with respect to the measure $\mu$.  Let $f: (0, \infty) \rightarrow  \mathbb{R}$ be a convex function.
The $*$-adjoint function $f^*:(0, \infty) \rightarrow  \mathbb{R}$ of $f$  is defined by (e.g. \cite{LieseVajda2006})
\begin{equation}\label{adjoint}
f^*(t) = t f (1/t), \ \  t\in(0, \infty).
\end{equation}
\par
It is obvious   that $(f^*)^*=f$ and that $f^*$ is again convex if $f$ is convex.
Csisz\'ar \cite{Csiszar}, and independently Morimoto \cite{Morimoto1963} and Ali \& Silvery \cite{AliSilvery1966} introduced 
the $f$-divergence   $D_f(P,Q)$ of the measures $P$ and $Q$ which, for a convex function $f: (0, \infty) \rightarrow  \mathbb{R}$ can be  defined as  (see \cite{LieseVajda2006})
\begin{eqnarray}\label{def:fdiv1}
D_f(P,Q)&=&
 \int_{\{pq>0\} }f\left(\frac{p}{q} \right) q d\mu + f(0)\  Q\left(\{x\in X: p(x) =0\}\right)\nonumber \\
 &+& f^*(0) \ P\left(\{x\in X: q(x) =0\}\right),
\end{eqnarray}
where 
\begin{equation}\label{fat0}
f(0) = \lim_{t\downarrow 0} f(t)  \  \  \text{ and} \  \   f^*(0) = \lim_{t\downarrow 0} f^*(t).
\end{equation}
We make the convention that $0 \cdot \infty =0$. 
\vskip 2mm
Please note that 
\begin{equation}\label{fstern}
D_f(P,Q)=D_{f^*}(Q,P).
\end{equation}
With (\ref{fat0}) and as  
$$f^*(0) \ P\left(\{x\in X: q(x) =0\}\right) = \int _{\{q=0\}} f^*\left(\frac{q}{p} \right) p  d\mu =  \int _{\{q=0\}} f\left(\frac{p}{q} \right) q  d\mu,$$
we can write in short
\begin{equation}\label{def:fdiv2}
D_f(P,Q)=
 \int_{X} f\left(\frac{p}{q} \right) q d\mu.
\end{equation}
\vskip 4mm 
For particular choices of $f$ we get many common divergences. E.g. for $f(t) = t \ln t$ with  $*$-adjoint function $f^*(t) = - \ln t$, the $f$-divergence is  the classical {\em information divergence}, also called {\em Kullback-Leibler divergence} or {\em relative entropy} from $P$ to $Q$ (see \cite{CT2006})
\begin{equation}\label{relent}
 D_{KL}(P\|Q)= \int_{X} p \ln \frac{p}{q} d\mu.
\end{equation}
For   the convex or concave functions  $f(t) = t^\alpha$ we obtain the {\em Hellinger integrals} (e.g. \cite{LieseVajda2006})
\begin{equation}\label{Hellinger}
H_\alpha (P,Q) = \int _X  p^\alpha q^{1-\alpha} d\mu.
\end{equation}
Those are related to the 
R\'enyi divergence of order $\alpha$, $\alpha \neq 1$,  introduced by  R\'enyi \cite{Ren} (for $\alpha >0$) as 
\begin{equation}\label{renyi}
D_\alpha(P\|Q)=
\frac{1}{\alpha -1} \ln \left( \int_X p^\alpha q^{1-\alpha} d\mu \right)= \frac{1}{\alpha -1} \ln \left( H_\alpha (P,Q)\right).
\end{equation}
The case $\alpha =1$ is the relative entropy $ D_{KL}(P\|Q)$.

\vskip 5mm
\section{$f$-divergences for  convex bodies.}

We will now  consider  $f$-divergences  for  convex bodies $K \in \mathcal{K}_0$. Let 
\begin{equation}\label{densities}
p_K(x)= \frac{ \kappa_{K}(x)}{\langle x, N_{K}(x) \rangle^{n} \  n|K^{\circ}|} \, , \   \ q_K(x)= \frac{\langle x, N_{K}(x) \rangle }{n\  |K|}.
\end{equation}
\noindent 
Usually, in the literature, 
the measures  under consideration are probability measures. Therefore
we have normalized the densities.
Thus 
\begin{equation}\label{PQ}
P_K=p_K\  \mu_K \ \ \ \text{and}   \ \ \   Q_K=q_K \ \mu_K
\end{equation}
are  measures on $\partial K$ that are absolutely continuous with respect  to $\mu_{K}$. $Q_K$ is a probability measure
and $P_K$ is one if $K$ is in $C^2_+$.
\par
Recall  that the normalized cone measure $cm_K$
on $\partial K$ is defined as follows:
For every measurable set $A \subseteq \partial K$
\begin{equation}\label{def:conemeas} 
cm_{K}(A)  = \frac{1}{|K|} \bigg| \big\{ta : \ a \in A, t\in [0,1] \big\} \bigg|.
\end{equation}
\vskip 3mm
\noindent
The next proposition is well known. See e.g. \cite {PaourisWerner2011} for a proof. It shows that the measures $P_K$ and $Q_K$ defined in (\ref{PQ})
 are the cone measures  of $K$  and $K^\circ$.  $N_K:\partial K \rightarrow S^{n-1}$, $x \rightarrow N_K(x)$  is the Gauss map.
\vskip 3mm
\noindent 
\begin{proposition} \label{prop:conemeas}
\noindent 
Let $K$ be a  convex body in $\mathbb R^n$. Let $P_K$ and $Q_K$ be the probability measures on $\partial K$  defined by (\ref{PQ}).  Then
$$
Q_K= cm_{K},
$$
\noindent 
or, equivalently, for every measurable subset $A$ in $ \partial K$
$Q_K(A)= cm_{ K}(A)$.
\newline
If $K$ is in addition in $C^2_+$, then 
$$
P_K= N_{K}^{-1}N_{K^{\circ}}cm_{K^{\circ}}
$$
\noindent 
or, equivalently, for every measurable subset $A$ in $ \partial K$
\begin{equation}\label{cone}
P_K(A)= cm_{ K^{\circ}} \bigg(N_{{K^{\circ}}}^{-1} \big(N_{K} (A)\big)\bigg).
\end{equation}
\end{proposition}
It is in the sense (\ref{cone}) that we understand $P_K$ to be the ``cone measure" of $K^\circ$ and we write $P_K=cm_{K^\circ}$.
\vskip 4mm 
 We now define the $f$-divergences of $K \in \mathcal{K}_0$. Note that
 $\langle x, N_K(x) \rangle > 0 $ for all $x \in \partial K$ and therefore $\{x \in \partial K: q_K(x) =0 \} = \emptyset$.
 Hence, possibly also using our convention  $0 \cdot \infty =0$, 
 $$f^*(0) \  P_K\left( \{x \in \partial K: q_K(x) =0 \} \right) = 0.$$

 \begin{definition}\label{f-div} 
{\em  Let $K$ be a convex body in $\mathcal{K}_0$ and let Let $f: (0, \infty) \rightarrow \mathbb{R}$ be a convex function. The 
$f$-divergence of $K$ with respect to the cone measures $P_K$ and $Q_K$ is
\begin{eqnarray}\label{f-div1,0}
D_f(P_K, Q_K)&= & \int _{\partial K} f\left(\frac{p_K}{q_K}\right) q_K d\mu_K \nonumber \\
&=& \int _{\partial K} f\left(\frac{|K| \kappa_K(x)} {|K^\circ| \langle x, N_K(x) \rangle ^{n+1}}\right)\frac{ \langle x, N_K(x) \rangle}{n |K|
} d\mu_K.
\end{eqnarray}}
\end{definition}
\vskip 2mm
\noindent
\vskip 4mm
\noindent
{\bf Remarks.}
\vskip 2mm
By  (\ref{fstern})  and (\ref{f-div1,0}) 
\begin{eqnarray}\label{f-div2,0}
 D_f(Q_K,P_K)&=& \int _{\partial K} f\left(\frac{q_K}{p_K}\right) p_K d\mu_K
= D_{f^*}(P_K, Q_K)  \nonumber \\
&= &\int _{\partial K} f^*\left(\frac{p_K}{q_K}\right) q_K d\mu_K  \nonumber \\
&=& \int _{\partial K} f\left(\frac{ |K^\circ| \langle x, N_K(x) \rangle ^{n+1}}{|K| \kappa_K(x)} \right)\frac{\kappa_K(x)\  d\mu_K} {n |K^\circ| \langle x, N_K(x) \rangle ^{n}}. 
\end{eqnarray}
\vskip 2mm
$f$-divergences can also be expressed as integrals over $S^{n-1}$, 
\begin{eqnarray}\label{f-div1,2}
D_f(P_K, Q_K)&= &
\int _{S^{n-1}} f\left(\frac{|K| } {|K^\circ| f_K(u) h_K(u) ^{n+1}}\right)\frac{ h_K(u)f_K(u)}{n |K|} d\sigma 
\end{eqnarray}
and
\vskip 2mm
\noindent
\begin{eqnarray}\label{f-div2,2}
D_f(Q_K,P_K)&=& 
\int _{ S^{n-1}} f\left(\frac{ |K^\circ| f_K(u) h_K(u) ^{n+1}}{|K| } \right)\frac{ d\sigma_K} {n |K^\circ| h_K(u)^{n}}.
\end{eqnarray}
\vskip 4mm
\noindent
{\bf Examples.}
\vskip 2mm
If $K$ is a polytope, the Gauss curvature $\kappa_K$ of  $K$ is $0$ a.e. on $\partial K$. Hence   
\begin{equation}\label{polytope}
D_f(P_K,Q_K)= f(0)  \  \   \text{and} \  \  D_f(Q_K,P_K) = f^*(0).
\end{equation}
\vskip 2mm
For every ellipsoid $\mathcal{E}$,  
\begin{equation}\label{ellipsoid}
D_f(P_{\mathcal{E}}, Q_{\mathcal{E}}) = D_f(Q_{\mathcal{E}}, P_{\mathcal{E}})=f(1) = f^*(1).
\end{equation}
\vskip 2mm
Denote by $Conv(0, \infty)$ the set of functions $\psi: (0, \infty) \rightarrow (0, \infty)$ such that $\psi$ is convex, $\lim_{t\rightarrow 0} \psi(t) = \infty$, and $\lim_{t\rightarrow \infty} \psi(t) = 0$. For $\psi \in Conv(0, \infty)$, Ludwig   \cite{Ludwig2010} 
introduces  the $L_\psi$ affine surface area for a convex body $K$ in $\mathbb{R}^n$ 
\begin{equation}\label{Lpsi}
\Omega_\Psi(K)= \int_{\partial K} \psi\left(\frac{ \kappa_K(x)} { \langle x, N_K(x) \rangle ^{n+1}}\right) \langle x, N_K(x) \rangle d \mu_K.
\end{equation}
Thus,  $L_\psi$ affine surface areas are special cases of   (non-normalized) $f$-divergences for $f=\psi$.
\vskip 2mm
For $\psi \in Conv(0, \infty)$, the $*$-adjoint function  $\psi^*$ is convex, $\lim_{t\rightarrow 0} \psi(t) = 0$, and $\lim_{t\rightarrow \infty} \psi(t) = \infty$.
Thus $\psi^*$ is an {\em Orlicz function} (see \cite{LindenstraussTzafriri}), and gives  rise to the corresponding {\em Orlicz-divergences}
$D_{\psi^*}(P_K,Q_K)$ and  $D_{\psi^*}(Q_K,P_K)$.
\vskip 2mm
Let $p \leq 0$. Then the function $f: (0, \infty) \rightarrow (0, \infty)$, $f(t)=t^\frac{p}{n+p}$,  is convex. The corresponding (non-normalized)
$f$-divergence (which is also an Orlicz-divergence) is the {\em $L_p$ affine surface area}, introduced by Lutwak \cite{Lutwak1996} for $p>1$
and by Sch\"utt and Werner \cite{SW2004} for $p< 1, p \neq -n$.  See also \cite{Hug}.
\par
It was shown in \cite{Werner2012} that all $L_p$ affine surface areas
are entropy powers of  R\'enyi divergences.
\vskip 2mm
For $p \geq 0$, the function $f: (0, \infty) \rightarrow (0, \infty)$, $f(t)=t^\frac{p}{n+p}$ is concave. The corresponding 
$L_p$ affine surface areas $\int_{\partial K} \frac{\kappa_K^\frac{p}{n+p}d \mu_K}{\langle x, N_K(x)\rangle ^\frac{n(p-1)}{n+p} } $
are examples of $L_\phi$ affine surface areas which were considered in \cite{LR2} and \cite{Ludwig2010}. Those, in turn are special cases
of (non-normalized) $f$-divergences for {\em concave} functions $f$.
\vskip 2mm
Let $f(t)= t \ln t$. Then the  $*$-adjoint function is $f^*(t) =- \ln t$. The corresponding $f$-divergence is the {\em Kullback Leibler divergence} or {\em relative entropy} $D_{KL}(P_K\|Q_K)$ from $P_K$ to $Q_K$
\begin{eqnarray}\label{KL1}
D_{KL}(P_K\|Q_K)  = \int_{\partial K}
\frac{\kappa_K(x)} {n |K^\circ| \langle x, N_{K}(x) \rangle ^n} \ln \left(\frac{|K|\kappa_K(x)}{|K^\circ|\langle x, N_{K}(x) \rangle ^{n+1}}\right) d \mu_K.
\end{eqnarray}
The relative entropy $D_{KL}(Q_K\|P_K)$ from $Q_K$ to $P_K$ is
\begin{eqnarray}\label{KL2}
D_{KL}(Q_K\|P_K) &= &D_{f^*}(P_K, Q_K) \\
&=&\int_{\partial K} \frac{\langle x, N_{K}(x) \rangle}{n|K|} \log\left(\frac{|K^\circ|\langle x, N_{K}(x) \rangle ^{n+1}}{|K| \kappa_K(x)} \right)d \mu_K.
\end{eqnarray}
Those were studied in detail in \cite{PaourisWerner2011}.
\vskip 4mm
Equations (\ref{f-div1,2}) and (\ref{f-div2,2}) of the above  remark   lead us to define $f$-divergences for several convex bodies, or mixed $f$-divergences.
\par
Let $K_1, \dots , K_n$ be convex bodies in $\mathcal{K}_0 $. 
Let $u \in S^{n-1}$. For $1 \leq i \leq n$, define
\begin{equation}\label{i-densities}
p_{K_i}(u)= \frac{1}{ n |K_i^{\circ}| h_{K_i}(u)} \, , \   \ q_{K_i}(u)= \frac{f_{K_i}(u) h_{K_i}(u) }{n\  |K_i|}.
\end{equation}
and measures on $S^{n-1}$ by
\begin{equation}\label{iPQ}
P_{K_i}=p_{K_i}\  \sigma \ \ \ \text{and}   \ \ \   Q_{K_i}=q_{K_i} \ \sigma.
\end{equation}
\vskip 4mm
\noindent
Let $f_i: (0, \infty) \rightarrow \mathbb{R}$, $ 1 \leq i \leq n$, be convex functions.
Then we  define the {\em mixed $f$-divergences}  for  convex bodies $K_1, \dots,  K_n$  in $\mathcal{K}_0$ by
\vskip 4mm
\noindent
\begin{definition}\label{mixed}
{\em 
\begin{equation*}
D_{f_1 \dots f_n}(P_{K_1} \times \dots \times P_{K_n}, Q_{K_1} \times \dots \times Q_{K_n})=  \int_{S^{n-1}} \prod_{i=1}^n 
\left[f_i\left(\frac{p_{K_i}}{q_{K_i}}\right) q_{K_i}\right]^\frac{1}{n}  d\sigma
\end{equation*}
and 
\begin{equation*}
D_{f_1\dots f_n}(Q_{K_1} \times \dots \times Q_{K_n}, P_{K_1} \times \dots \times P_{K_n})=  \int_{S^{n-1}} \prod_{i=1}^n \left[f_i\left(\frac{q_{K_i}}{p_{K_i}}\right) p_{K_i} \right]^\frac{1}{n} d\sigma.
\end{equation*}}
\end{definition}
Note that  
\begin{eqnarray*}
&&\hskip -15mm D_{f_1^* \dots f_n^*}(P_{K_1} \times \dots \times P_{K_n}, Q_{K_1} \times \dots \times Q_{K_n})\\
&& \hskip 20mm= 
D_{f_1\dots f_n}(Q_{K_1} \times \dots \times Q_{K_n}, P_{K_1} \times \dots \times P_{K_n}).
\end{eqnarray*}
\vskip 4mm
Here, we concentrate on $f$-divergence for one convex body. Mixed $f$-divergences are  treated similarly.
We also refer to \cite{Ye2012}, where they have been investigated for functions in $Conv(0,\infty)$.

\vskip 4mm
The observation (\ref{polytope}) about polytopes  holds more generally.
\vskip 3mm
\begin{proposition}\label{prop1}
Let $K$ be a convex body in $\mathcal{K}_0$  and let  $f: (0, \infty) \rightarrow \mathbb{R}$ be a convex function. 
If $K$ is such that $\mu_K\left(\{p_K>0\}\right) =0$, then 
$$
D_f(P_K,Q_K) = f(0)  \ \  \text{and} \  \  D_f(Q_K, P_K) = f^*(0).
$$ 
\end{proposition}
\vskip 2mm
\noindent
{\bf Proof.}
$\mu_K\left(\{p_K>0\}\right) =0$ iff $Q_K \left(\{p_K>0\}\right) =0$. Hence  the assumption 
implies that $Q_K \left(\{p_K=0\}\right) =1$. Therefore, 
\begin{eqnarray*}
D_f(P_K, Q_K)&= &
\int _{\partial K} f\left(\frac{p_K}{q_K}\right) q_K d\mu_K \\
&=& \int _{\{p_K>0\} } f\left(\frac{p_K}{q_K}\right) q_K d\mu_K + \int _{\{p_K=0\} } f\left(\frac{p_K}{q_K}\right) q_K d\mu_K\\
&=& f(0).
\end{eqnarray*}
By (\ref{fstern}),  $D_f(Q_K, P_K)= D_{f^*}(P_K, Q_K)=f^*(0)$.

\vskip 3mm
The next proposition complements the previous one. 
In view of  (\ref{ellipsoid}) and (\ref{affin3}), it corresponds to the affine isoperimetric inequality for $f$-divergences. It was proved in \cite{LieseVajda2006} in a different setting 
and in the special case of $f \in Conv(0, \infty)$ by Ludwig \cite{Ludwig2010}.
We include a  proof for completeness.
\vskip 3mm
\begin{proposition}\label{prop1}
Let $K$ be a convex body in $\mathcal{K}_0$  and let  $f: (0, \infty) \rightarrow \mathbb{R}$ be a convex function. 
If $K$ is such that $\mu_K\left(\{p_K>0\}\right) >0$, then 
\begin{equation}\label{affin1}
D_f(P_K,Q_K) \geq f \left(\frac{P_K\left(\{p_K>0\} \right)}{Q_K \left(\{p_K>0\}\right)}\right) \  Q_K \left(\{p_K>0\}\right) + f(0)  \  Q_K \left(\{p_K=0\}\right)
\end{equation}
and
\begin{equation}\label{affin2}
D_f(Q_K, P_K) \geq  f^* \left(\frac{P_K\left(\{p_K>0\} \right)}{Q_K \left(\{p_K>0\}\right)}\right) \  Q_K \left(\{p_K>0\}\right) + f^*(0) \  Q_K \left(\{p_K=0\}\right).
\end{equation}
If $K$ is  in $C^2_+$,  or if $f$ is decreasing, then 
\begin{equation}\label{affin3}
D_f(P_K, Q_K) \geq  f(1)  \   \   \text{and}  \   \   D_f(Q_K, P_K) \geq  f^*(1)=f(1).
\end{equation}
Equality holds  in (\ref{affin1}) and   (\ref{affin2})  iff  $f$ is linear or $K$ is an ellipsoid.
If $K$ is in $C^2_+$, equality holds  in both inequalities (\ref{affin3}) iff   $f$ is linear or $K$ is an ellipsoid. If  $f$ is decreasing, 
equality holds  in  both inequalities (\ref{affin3}) iff $K$ is an ellipsoid.
\end{proposition}
\vskip 2mm
\noindent
{\bf Remark.} It is possible for  $f$ to be deceasing and linear without having equality in  (\ref{affin3}). 
To see that, let $f(t) =at +b$, $a < 0$, $b>0$. Then, for polytopes $K$ (for which $\mu_K\left(\{p_K>0\}\right)=0$), $D_f(P_K, Q_K) = f(0) = b > f(1) = a+b$. 
But, also in the case when $0 < \mu_K\left(\{p_K>0\}\right) <1$, strict inequality may hold. 
\par
Indeed, let $\varepsilon >0$ be sufficiently small and let $K= B^n_\infty(\varepsilon)$ be a ``rounded" cube, 
where we have ``rounded" the corners of the cube $B^n_\infty$ with sidelength  $2$ centered at $0$ by replacing each corner with $\varepsilon B^n_2$ Euclidean balls. Then  $D_f(P_K, Q_K) = b  + a \  P_K\left(\{p_K>0\} \right) > b + a = f(1)$.
\vskip 3mm
\noindent
{\bf Proof of Proposition \ref{prop1}.}
Let $K$ be such that $\mu_K\left(\{p_K>0\}\right) >0$, which is equivalent to  $Q_K \left(\{p_K>0\}\right) >0$. Then, 
by Jensen's inequality, 
\begin{eqnarray*}
D_f(P_K, Q_K)&=&  Q_K \left(\{p_K>0\}\right) \ \int _{\{p_K>0\} } f\left(\frac{p_K}{q_K}\right) \frac{q_K d\mu_K}{Q_K \left(\{p_K>0\}\right)} \\
&+ &f(0) \  Q_K \left(\{p_K=0\}\right) \\
& \geq& Q_K \left(\{p_K>0\}\right)   f\left( \frac{P_K \left(\{p_K>0\}\right)}{Q_K \left(\{p_K>0\}\right)} \right) + f(0) Q_K \left(\{p_K=0\}\right).
\end{eqnarray*}
Inequality (\ref{affin2}) follows  by (\ref{fstern}), as
$
D_f(Q_K, P_K) = D_{f^*} (P_K, Q_K)$.
\vskip 2mm
If $K$ is in $C^2_+$, $Q_K \left(\{p_K>0\}\right)=1$, $Q_K \left(\{p_K=0\}\right)=0$, $P_K \left(\{p_K>0\}\right)=1$ and $P_K \left(\{p_K=0\}\right)=0$. Thus we get that $D_f(P_K, Q_K) \geq f(1)$ and $D_f(Q_K, P_K) \geq f^*(1)=f(1)$.
\par
If $f$ is decreasing, then, by  Jensen's inequality   
\begin{eqnarray*}
D_f(P_K, Q_K)&=&   \int _{\partial K } f\left(\frac{p_K}{q_K}\right)q_K d\mu_K \geq f \left(  \int _{\partial K }p_K d\mu_K\right) \geq f(1).
\end{eqnarray*}
The last inequality holds as $ \int _{\partial K }p_K d\mu_K \leq 1$ and as $f$ is decreasing.
\vskip 2mm
Equality holds in Jensen's inequality iff either $f$ is linear or $\frac{p_K}{q_K}$ is constant.
Indeed, if $f(t) = at +b$, then 
\begin{eqnarray*}
D_f(P_K, Q_K)&=& 
\int _{\{p_K>0\} } \left(a \frac{p_K}{q_K} +b \right) q_K d\mu_K +  f(0) \  Q_K \left(\{p_K=0\}\right)\\
&=& a P_K \left(\{p_K>0\}\right) + f(0).
\end{eqnarray*}
If $f$ is not linear, equality holds iff $\frac{p_K}{q_K}=c$,  $c$ a  constant. As by assumption $\mu_K\left(\{p_K>0\}\right) >0$, 
$c \neq 0$. By a theorem of Petty \cite{Petty}, this holds iff $K$ is an ellipsoid.
\vskip 4mm
The next proposition can be found  in \cite{LieseVajda2006} in a different setting.  Again, we include a proof for completeness.
\begin{proposition}\label{prop2}
Let $K$ be a convex body in $\mathcal{K}_0$  and let  $f: (0, \infty) \rightarrow \mathbb{R}$ be a convex function. Then
$$
D_f(P_K,Q_K) \leq f(0) + f^*(0) + f(1)\bigg[ Q_K(\{0 <p_K \leq q_K\} ) + P_K(\{0 < q_K\leq p_K\})\bigg]
$$
and 
$$
D_f(Q_K, P_K) \leq f(0) + f^*(0) + f(1)\bigg[ Q_K(\{0 <p_K \leq q_K\} ) + P_K(\{0 < q_K\leq p_K\})\bigg].
$$
If $f$ is decreasing, the inequalities reduce to $D_f(P_K, Q_K) \leq f(0)$ respectively, $D_f(Q_K, P_K) \leq f^*(0)$.
\end{proposition}
\vskip 2mm
\noindent
{\bf Proof.} It is enough to prove the first inequality. The second one follows immediately form the first  by (\ref{fstern}).

\begin{eqnarray*}
D_f(P_K, Q_K)&= &
\int _{\partial K} f\left(\frac{p_K}{q_K}\right) q_K d\mu \\
&=&  \int _{\{p_K>0\}} f\left(\frac{p_K}{q_K}\right) q_K d\mu + f(0) \  Q_K(\{p_K=0\})\\
&=& f(0) \  Q_K(\{p_K=0\}) +  \int _{\{0 < p_K\}  \cap \{f^\prime \geq 0\}} f\left(\frac{p_K}{q_K}\right) q_K d\mu \\
&&\hskip 35mm +
  \int _{\{0 < p_K\}  \cap \{f^\prime \leq 0\}} f\left(\frac{p_K}{q_K}\right) q_K d\mu \\
&\leq &f(0) \  \bigg[Q_K(\{p_K=0\}) + Q_K\left(\{p_K>0\} \cap \{f^\prime \leq 0\} \right) \bigg]   \\
&&\hskip -25mm
+ \int _{\{0 < p_K\leq q_K \}  \cap \{f^\prime \geq 0\}} f\left(\frac{p_K}{q_K}\right) q_K d\mu 
+
 \int _{\{0 < q_K\leq p_K\} \cap \{f^\prime \geq 0\} } f\left(\frac{p_K}{q_K}\right) q_K d\mu \\
&\leq&  f(0) + f(1) \ Q_K\left(\{0 <p_K \leq q_K\}  \cap \{f^\prime \geq 0\}\right)   \\ 
&+& \int _{\{0 < q_K\leq p_K\}  \cap \{f^\prime \geq 0\}} f^*\left(\frac{q_K}{p_K}\right) p_K d\mu \\
&=&  f(0) + f(1)  \ Q_K\left(\{0 <p_K \leq q_K\}  \cap \{f^\prime \geq 0\}\right)  \\
& +&
\int _{\{0 < q_K\leq p_K\}  \cap \{f^\prime \geq 0\} \cap \{(f^*)^\prime \geq 0\} } f^*\left(\frac{q_K}{p_K}\right) p_K d\mu \\
& +&
\int _{\{0 < q_K\leq p_K\}  \cap \{f^\prime \geq 0\} \cap \{(f^*)^\prime \leq 0\} } f^*\left(\frac{q_K}{p_K}\right) p_K d\mu\\
& \leq&  f(0) + f(1)  \ Q_K\left(\{0 <p_K \leq q_K\}  \cap \{f^\prime \geq 0\}\right)  \\
& +&  f^*(1)  \ P_K\left(\{0 <q_K \leq p_K\}  \cap \{f^\prime \geq 0\} \cap \{(f^*)^\prime \geq 0\} \right)  \\
& + & f^*(0) \ P_K\left(\{0 <q_K \leq p_K\}  \cap \{f^\prime \geq 0\} \cap \{(f^*)^\prime \leq 0\} \right)  \\
&&\hskip -29mm \leq f(0)  + 
  f^*(0) \ P_K\left(\{0 <q_K \leq p_K\}  \cap \{f^\prime \geq 0\} \right) \\
&&\hskip -29mm + f(1)\bigg[ Q_K(\{0 <p_K \leq q_K\}\cap \{f^\prime \geq 0\} ) + P_K(\{0 < q_K\leq p_K\}\cap \{f^\prime \geq 0\})\bigg].
\end{eqnarray*}
It follows from the last expression that, if $f$ is decreasing, the inequality reduces to $D_f(P_K, Q_K) \leq f(0)$.
\vskip 4mm
The next  proposition shows that $f$-divergences are $GL(n)$ invariant and that non-normalized  
$f$-divergences are $SL(n)$ invariant valuations.
For functions in $Conv(0, \infty)$, this was proved by Ludwig \cite{Ludwig2010}. 
\par
For functions in $Conv(0, \infty)$ the expressions are also lower semicontinuous,  as it was shown in 
\cite{Ludwig2010}.  However, this  need not be the case anymore
if we assume just convexity of $f$. Indeed, let $f(t)=t^2$ and let $K=B^n_2$ be the Euclidean unit ball.
Let $(K_j)_{j \in \mathbb{N}}$ be a sequence of polytopes that converges to $B^n_2$.
As observed above, $D_f(P_{K_j},Q_{K_j})=f(0) =0$ for all $j$. But $D_f(P_{B^n_2},Q_{B^n_2})=f(1) =1$.
\par 
Let $\tilde{P}_K = \frac{\kappa_K   \mu_K}{\langle x, N_K(x)\rangle ^n} $ and $\tilde{Q}_K = \langle x, N_K(x)\rangle  \mu_K$.
Then we will denote by $D_f(\tilde{P}_K, \tilde{Q}_K)$ and $D_f(\tilde{Q}_K, \tilde{P}_K)$ the non-normalized $f$-divergences.
We will also use the following lemma from \cite{SW2004} for the proof of Proposition \ref{prop3}.
\vskip 3mm
\begin{lemma}\label{lem:transform}
Let $K$ be a convex body in $\mathcal{K}_0$. Let  $h:\partial K\rightarrow
\mathbb R$ be an integrable function, and $T:\mathbb R^{n}\rightarrow
\mathbb R^{n}$
an invertible, linear map. Then
$$
\int_{\partial K}h(x)d\mu_{K}
=|\det(T)|^{-1}
\int_{\partial T(K)}
\frac{f(T^{-1}(y))}{\|T^{-1t}(N_{K}(T^{-1}(y)))\|}d\mu_{T(K)}.
$$
\end{lemma}

\vskip 3mm
\begin{proposition}\label{prop3}
Let $K$ be a convex body in $\mathcal{K}_0$ and let  $f: (0, \infty) \rightarrow \mathbb{R}$ be a convex function. Then $D_f(P_K,Q_K)$ and $D_f(Q_K, P_K)$ are $GL(n)$ invariant and $D_f(\tilde{P}_K, \tilde{Q}_K)$ and $D_f(\tilde{Q}_K, \tilde{P}_K) $ are $SL(n)$ invariant valuations.
\end{proposition}
\vskip 2mm
\noindent
{\bf Proof.} 
We use (e.g. \cite{SW2004}) that
\begin{eqnarray*}
\langle T(x), N_{ T(K)}(T(x) \rangle
=\frac{\langle x,N_{ K}(x)\rangle }
{\|T^{-1t}(N_{ K}(x))\|},
\end{eqnarray*}
and
$$
\kappa_K(x)^{}
=\|T^{-1t}(N_{ K}(x))\|^{n+1}\det(T)^{2}
\kappa_{T(K)}(T(x))^{}
$$
and Lemma \ref{lem:transform} to get that
\begin{eqnarray*}
D_f(P_K,Q_K)&=& \int_{\partial K}f\left(\frac{p_K(x)}{q_K(x)} \right) q_K(x) d \mu(x) \\
&=& \frac{1}{|\det(T)|}
\int_{\partial T(K)}\frac
{f\left(\frac{p_{K}(T^{-1}(y))}{q_{K}(T^{-1}(y))}\right) q_{K}(T^{-1}(y)) d\mu_{T(K)}}{\|T^{-1t}(N_{ K}(T^{-1}(y)))\|} \\
&=& D_f(P_{T(K)},Q_{T(K)}).
\end{eqnarray*}
The formula for  $D_f(Q_K,P_K)$ follows immediately from this one and  (\ref{fstern}).
The $SL(n)$ invariance for the non-normalized $f$-divergences is shown in the same way.
\par
Now we show that $D_f(\tilde{P}_K,\tilde{Q}_K)$ and $D_f(\tilde{Q}_K, \tilde{P}_K)$ are valuations, i.e.
for convex bodies $K$ and $L$ in $\mathcal{K}_0$ such
that $K \cup L \in \mathcal{K}_0$,
\begin{equation}\label{val}
D_f(\tilde{P}_{K\cup L},\tilde{Q}_{K\cup L}) + D_f(\tilde{P}_{K\cap L},\tilde{Q}_{K\cap L}) = D_f(\tilde{P}_K,\tilde{Q}_K) + D_f(\tilde{P}_L,\tilde{Q}_L).
\end{equation}
Again, it is enough to prove this formula and the one for $D_f(\tilde{Q}_K,\tilde{P}_K)$ follows with  (\ref{fstern}).
To prove (\ref{val}), we proceed as in
Sch\"utt  \cite{Schuett1993}. For completeness, we  include the argument. We decompose
$$
\partial (K \cup L) = (\partial K \cap \partial L)  \cup (\partial K \cap  L^c) \cup ( K^c \cap \partial L),
$$
$$
\partial (K \cap L) = (\partial K \cap \partial L)  \cup (\partial K \cap \text{int} L) \cup (\text{int} K \cap \partial L),
$$
$$
\partial K = (\partial K \cap \partial L)  \cup (\partial K \cap L^c) \cup (\partial K \cap \text{int} L),
$$
$$
\partial L = (\partial K \cap \partial L)  \cup (\partial K^c \cap \partial L) \cup (\text{int} K \cap \partial L),
$$
where all unions on the right hand side are disjoint. Note that for $x$ such
that the curvatures $\kappa_K(x)$, 
$\kappa_L(x)$, $\kappa_{K\cup L} (x)$ and $\kappa_{K\cap L} (x)$ exist,
\begin{equation}\label{support}
\langle x, N_K(x)\rangle  = \langle x, N_L(x) \rangle =  \langle x, N_{K\cap L}(x) \rangle = \langle x,  N_{K\cup L}(x) \rangle
\end{equation}
and
\begin{equation}\label{curv}
\kappa_{K\cup L}(x)= \min\{\kappa_K(x), \kappa_L(x)\}, \  \  \kappa_{K\cap L} (x)= \max\{\kappa_K(x), \kappa_L(x)\}.
\end{equation}
To prove (\ref{val}), we split the involved integral using the above decompositions  
and  (\ref{support}) and (\ref{curv}).
\vskip 5mm
\section{Geometric characterization of $f$-divergences.}
In \cite{Werner2012}, geometric characterizations were proved  for  R\'enyi divergences. Now, we want to establish such  geometric characterizations for $f$-divergences as well.
We use the {\em surface body} \cite{SW2004}  but the
{\em illumination surface body} \cite{WernerYe2010} or the {\em mean width body}  \cite{JenkinsonWerner2012} can also be used.
 \vskip 5mm
Let $K$ be a convex body in $\mathbb R^{n}$.
Let $g:\partial K\rightarrow \Bbb R$  be a nonnegative,  integrable, function. 
Let $s\geq 0$. 
\par
The {\em surface body}  $K_{g,s}$,  introduced in \cite{SW2004},  is the intersection of all closed half-spaces
$H^{+}$
whose defining hyperplanes $H$ cut off a set of 
$f \mu_K$-measure less
than or equal to $s$ from $\partial K$. More precisely,
\begin{equation*}\label{def1}
K_{g,s}= \bigcap _{\int_{\partial K \cap H^-}g d\mu_K \leq s} H^+.
\end{equation*}
\vskip 2mm
For $ x \in \partial K$ and $s > 0$ 
$$x_s=[0,x] \cap \partial K_{g,s}.
$$

The {\em minimal function} $M_g : \partial K \rightarrow \mathbb{R}$
\begin{equation}\label{min}
M_{g}(x)=\inf_{0<s} \  \frac{\int_{\partial K\cap
H^{-}(x_{s},N_{ K_{g,s}}(x_{s}))}g\ d\mu_{ K}}{\mbox{vol}_{n-1}\left(\partial K\cap
H^{-}(x_{s},N_{ K_{g,s}}(x_{s}))\right)}
\end{equation}
was introduced  in   \cite{SW2004}.  
$H(x,\xi)$ is the hyperplane through $x$ and orthogonal to
$\xi$.
$H^{-}(x,\xi)$ is the closed halfspace containing the point $x+\xi$,
$H^{+}(x,\xi)$ the other halfspace.

For $x \in \partial K$, we define $r(x)$ as
the maximum of all real numbers $\rho$ so that
$B_{2}^{n}(x-\rho N_{ K}(x),\rho)\subseteq K$. 
Then we formulate an  integrability condition for the minimal function
\begin{equation}\label{intcond}
\int_{\partial K}\frac{d\mu_{
K}(x)}{\left(M_{g}(x)\right)^{\frac{2}{n-1}}r(x)}
<\infty.
\end{equation}

\vskip 4mm
The following theorem was proved in \cite{SW2004}.
\vskip 3mm
\begin{theorem}\label{thm:surface}  
Let $K$ be a convex body in $\mathbb R^{n}$.
Suppose that $f:\partial K\rightarrow \mathbb R$ is an integrable,
almost everywhere strictly positive function
that satisfies the integrability condition (\ref{intcond}).
Then
$$
c_{n}\lim_{s \to 0}
\frac{|K|-|K_{g,s}|}
{s^\frac{2}{n-1}}=
\int_{\partial K} \frac{\kappa_K^\frac{1}{n-1}}
{g^\frac{2}{n-1}}d\mu_{K},
$$
where $c_n=2 |B_2^{n-1}|^{\frac{2}{n-1}}$. 
\end{theorem}

\vskip 4mm
Theorem \ref{thm:surface}  was used  in \cite{SW2004}  to give  geometric interpretations of $L_p$ affine surface area
and in \cite{Werner2012} to give geometric interpretations of R\'enyi divergences.
Now we use this theorem to give 
geometric interpretations of $f$-divergence for cone measures of convex bodies.
\vskip 2mm
For a convex function $f:(0, \infty) \rightarrow \mathbb{R}$, 
 let $g_f, h_f: \partial K \rightarrow \mathbb{R}$  be defined as
\begin{equation}\label{gf}
g_f(x)= \left[n |K^\circ| n^n |K|^n  \frac{p_K q_K}{\left(f\left( \frac{p_K}{q_K}\right)\right)^{n-1}}\right]^\frac{1}{2}
\end{equation}
and
\begin{equation}\label{hf}
h_f(x)= g_{f^*}(x)= \left[n |K^\circ| n^n |K|^n  \frac{q_K^n/ p_K^{n-2}}{\left(f\left( \frac{p_K}{q_K}\right)\right)^{n-1}}\right]^\frac{1}{2}.
\end{equation}

\vskip 4mm
\begin{corollary}\label{cor:asp}
Let $K$ be a convex body in $\mathcal{K}_0$ and let $f:(0, \infty) \rightarrow \mathbb{R}$ be convex.
 Let $g_f, h_f: \partial K \rightarrow \mathbb{R}$  be defined as in (\ref{gf}) and (\ref{hf}).
 If  $g_f$  and $h_f$ are integrable,  almost everywhere strictly positive functions that satisfy the integrability condition (\ref {intcond}), then
 \begin{eqnarray*}
c_n \lim_{s \to 0}
\frac{|K|
-|K_{g_f,s}|}{  s^\frac{2}{n-1}}
= D_f(P_K,Q_K) 
\end{eqnarray*}
\vskip 2mm
\noindent
and 
\begin{eqnarray*}
c_n \lim_{s \to 0}
\frac{|K|
-|K_{h_f,s}|}{  s^\frac{2}{n-1}}
= D_f(Q_K,P_K) 
\end{eqnarray*}
\vskip 2mm
\noindent
\end{corollary}
\vskip 4mm
\noindent
{\bf Proof.} 
The proof of the corollary follows immediately from Theorem \ref{thm:surface}.

 \vskip 2mm \noindent Elisabeth Werner\\
{\small Department of Mathematics \ \ \ \ \ \ \ \ \ \ \ \ \ \ \ \ \ \ \ Universit\'{e} de Lille 1}\\
{\small Case Western Reserve University \ \ \ \ \ \ \ \ \ \ \ \ \ UFR de Math\'{e}matique }\\
{\small Cleveland, Ohio 44106, U. S. A. \ \ \ \ \ \ \ \ \ \ \ \ \ \ \ 59655 Villeneuve d'Ascq, France}\\
{\small \tt elisabeth.werner@case.edu}\\ \\

\end{document}